\documentclass[a4paper,11pt,oneside]{article}
\usepackage[english]{babel}
\usepackage{amsmath,amsfonts,latexsym,amssymb,amsthm}
\begin{document}
\begin{center}{\Large\bf Torsion of elliptic curves over quadratic cyclotomic fields}\\
\vspace{1cm}
Filip Najman
\end{center}
\vspace{1cm}
\begin{abstract} In this paper we study the possible torsions of elliptic curves over $\mathbb Q(i)$ and $\mathbb Q(\sqrt {-3})$.\\
\end{abstract}
\footnotetext{Mathematics subject classification (2000) 11G05, 14G05}
\textbf{1 Introduction}\\

For an elliptic curve $E$ over a number field $K$, it is well known, by the Mordell-Weil theorem, that the set $E(K)$ of $K$-rational points on $E$ is a finitely generated abelian group. The group $E(K)$ is isomorphic to $T\oplus\mathbb Z^r$, where $r$ is a non-negative integer and $T$ is the torsion subgroup. When $K=\mathbb Q$, by Mazur's Theorem, the torsion subgroup is either cyclic of order $m$, where $1 \leq m \leq 10$ or $m=12$, or of the form $\mathbb Z_2 \oplus \mathbb Z_{2m}$, where $1 \leq m \leq 4$.\\
If $K$ is a quadratic field, then the following theorem classifies the possible torsions.
\newtheorem*{tmx}{Theorem(Kamienny, \cite{kam}, Kenku and Momose, \cite{km})} 
\begin{tmx}
Let $K$ be a quadratic field and $E$ an elliptic curve over $K$. Then the torsion subgroup $E(K)_{tors}$ of $E(K)$ is isomorphic to one of the following $26$ groups:
$$\mathbb Z_m, \text{ for } 1 \leq m\leq 18,\ m\neq 17,$$
$$\mathbb Z_2 \oplus \mathbb Z_{2m}, \text{ for } 1 \leq m\leq 6,$$
$$\mathbb Z_3 \oplus \mathbb Z_{3m}, \text{ for }  m=1,2$$
$$\mathbb Z_4 \oplus \mathbb Z_{4}.$$
\end{tmx}
Moreover, the only quadratic field over which torsion $\mathbb Z_4 \oplus \mathbb Z_4$ appears is $\mathbb Q(i)$ and the only quadratic field over which torsions $\mathbb Z_3 \oplus \mathbb Z_3$ and $\mathbb Z_3 \oplus \mathbb Z_6$ appear is $\mathbb Q(\sqrt{-3})$.\\
In \cite{jkp}, Theorem 3.5, it is proved that if we let the quadratic fields vary, then all of the 26 torsion subgroups appear infinitely often.\\
In this paper we will take a different approach. We will fix the quadratic field and then study the possible torsions. The fields that we are going to study, $\mathbb Q(i)$ and $\mathbb Q(\sqrt {-3})$, are somewhat special, as they are the only fields containing roots of unity apart from $1$ and $-1$, i.e. they are the only cyclotomic quadratic fields. Also, over each of these fields, torsion subgroups appear that appear over no other fields. Note that the rings of integers of both these fields are unique factorization domains.\\
The results obtained for elliptic curves over $\mathbb Q(i)$ are presented in the following theorem.
\newtheorem{tm}{Theorem}
\begin{tm}
\begin{itemize}
\item[(i)] Let $E$ be an elliptic curve with rational coefficients. Then $E(\mathbb Q(i))_{tors}$ is either one of the groups from Mazur's Theorem or $\mathbb Z_4 \oplus \mathbb Z_4$.
\item[(ii)] Let $E$ be an elliptic curve defined over $\mathbb Q(i)$. Then $E(\mathbb Q(i))_{tors}$ is either one of the groups from Mazur's Theorem, $\mathbb Z_4 \oplus \mathbb Z_4$ or $\mathbb Z_{13}$.
\end{itemize}
\end{tm}
The first part of this theorem is the best possible, while for the second part we believe that $\mathbb Z_{13}$ does not appear as a torsion subgroup, but we were unable to prove this.\\
Note that the torsion subgroup $\mathbb Z_4 \oplus \mathbb Z_4$ appears infinitely often. Elliptic curves with this torsion can be written in the form $$y^2=x(x+m^2)(x+n^2),\ m,n\in\mathbb Z[i],$$ where $m^2-n^2$ is a square in $\mathbb Z[i]$. This is an easy corollary of the 2-descent proposition (see \cite{kn}, Theorem 4.2, p. 85).\\
The results obtained for elliptic curves over $\mathbb Q(\sqrt {-3})$ are presented in the following theorem.
\newtheorem{tm4}[tm]{Theorem}
\begin{tm4}
\begin{itemize}
\item[(i)] Let $E$ be an elliptic curve with rational coefficients. Then $E(\mathbb Q(\sqrt{-3}))_{tors}$ is either one of the groups from Mazur's Theorem, $\mathbb Z_3 \oplus \mathbb Z_3$ or $\mathbb Z_3 \oplus \mathbb Z_6$.
\item[(ii)] Let $E$ be an elliptic curve defined over $\mathbb Q(\sqrt{-3})$. Then $E(\mathbb Q(\sqrt{-3}))_{tors}$ is either one of the groups from Mazur's Theorem, $\mathbb Z_3 \oplus \mathbb Z_3$, $\mathbb Z_3 \oplus \mathbb Z_6$, $\mathbb Z_{13}$ or $\mathbb Z_{18}$.
\end{itemize}
\end{tm4}
Again, the first part of this theorem is the best possible($\mathbb Z_3 \oplus \mathbb Z_3$ and $\mathbb Z_3 \oplus \mathbb Z_6$ appear infinitely often), while for the second part we believe that $\mathbb Z_{13}$ and $\mathbb Z_{18}$ do not appear as torsion subgroups, but we were unable to prove this.\\
\vspace{0.3cm}\\
\textbf{2 Torsion over $\mathbb Q(i)$}\\

Throughout this chapter, the following extension of the Lutz-Nagell Theorem is used to compute torsion groups of elliptic curves. 
\newtheorem*{tm100}{Theorem (Extended Lutz-Nagell Theorem)}
\begin{tm100}
Let $E: y^2=x^3+Ax+B$ with $A,B\in \mathbb Z[i]$. If a point $(x,y)\in E(\mathbb Q(i))$ has finite order, then
\begin{enumerate}
\item $x,y\in \mathbb Z[i].$
\item Either $y=0$ or $y^2|4A^3+27B^2.$
\end{enumerate} 
\end{tm100}
The proof of the Lutz-Nagell Theorem can easily be extended to elliptic curves over $\mathbb Q(i)$. Details of the proof can be found in \cite{th}, Chapter 3. An implementation in Maple can be found in \cite{th}, Appendix A.\\
Next, we give a result that applies to elliptic curves over all quadratic fields and that is an immediate corollary of the main result of \cite{ll} (see also \cite{fu2}).
\newtheorem{tm3}[tm]{Lemma}
\begin{tm3}
Let $E$ be an elliptic curve with rational coefficients and $d$ a square-free integer. Then $E(\mathbb Q (\sqrt d))_{tors}$ cannot be $\mathbb Z_{11}, \mathbb Z_{13}$ or $\mathbb Z_{14}.$
\end{tm3}

Next we give a series of lemmas that will prove Theorem 1. 
\newtheorem{tm5}[tm]{Lemma}
\begin{tm5} $E(\mathbb Q(i))_{tors}$ cannot be $\mathbb Z_{11}, \mathbb Z_{18}$ or $\mathbb Z_{2}\oplus \mathbb Z_{10}.$
\end{tm5}
\emph{Proof:}\\
It is proved in \cite{km}, Example 3.2 that the torsion of an elliptic curve over $\mathbb Q(i)$ cannot be $\mathbb Z_2 \oplus \mathbb Z_{10}$.\\
It is proved in \cite{re}, Theorem 2 that the torsion cannot be $\mathbb Z_{11}$.\\
Since the rational prime 3 remains prime in $\mathbb Z[i]$, condition i) of Proposition 2.4 from \cite{km} is satisfied and hence, the torsion cannot be $\mathbb Z_{18}$. \qed

\newtheorem{tm6}[tm]{Lemma}
\begin{tm6} $E(\mathbb Q(i))_{tors}$ cannot be $\mathbb Z_{16}$.
\end{tm6}
\emph{Proof}:\\
By \cite{ra}, case 2.5.5., page 37, we see that elliptic curves with torsion $\mathbb Z_{16}$ over $\mathbb Q(i)$ are induced by the solutions of the equation \begin{equation}s^2=t(t^2+1)(t^2+2t-1),\ s,t\in \mathbb Q(i),\label{j1}\end{equation} where $t$ satisfies
\begin{equation}t(t^4-1)(t^2+2t-1)(t^2-2t-1)\neq 0.\label{j2}\end{equation}
It follows that it is enough to prove that this equation has no solutions. Let \begin{equation}t=\alpha\square,\label{eq2}\end{equation}
\begin{equation} t^2+1=\beta\square \label{eq3}\end{equation} and \begin{equation}t^2+2t-1=\gamma\square, \label{eq4}\end{equation} where $\alpha, \beta$ and $\gamma$ are square-free, nonzero Gaussian integers. Also, let $t=\frac u v$, where $u$ and $v$ are coprime, nonzero Gaussian integers.\\
First, we prove that $\alpha$ is relatively prime to $\beta$ and $\gamma$.
Suppose a Gaussian prime $\pi$ divides both $\alpha$ and $\beta$. From (\ref{eq2}) it follows that $\pi$ divides either $u$ or $v$ an odd number of times. From (\ref{eq3}), it follows that $\pi$ divides $u^2+v^2$, and since it divides exactly one of $u,v$, this is impossible.\\
Suppose a Gaussian prime $\pi$ divides both $\alpha$ and $\gamma$. Again, $\pi$ divides exactly one of $u,v$. From (\ref{eq4}), it follows that $\pi$ divides $u^2+2uv-v^2$, which is again impossible.\\
Since $-1$ is a square in $\mathbb Z[i]$, we conclude that $\alpha\in\{1,i\}$.\\
Next, we prove that $\gcd(\beta, \gamma)=1$ or $1+i$. Let $\pi$ be a Gaussian prime dividing both $\beta$ and $\gamma$. By subtracting (\ref{eq3}) from (\ref{eq4}), we conclude that $\pi$ divides $2uv-2v^2=2v(u-v)$. As it was already proved, since $\pi$ divides $\beta$, $\pi$ cannot divide $u$, implying $\pi|2(u-v)$. Suppose $\pi|(u-v)$, i.e. $u\equiv v \pmod \pi$.  Now, (\ref{eq3}) implies $2u^2\equiv 0\pmod \pi$, again implying $\pi|2$. Since $2=-i(1+i)^2$, we conclude that the only possibilities for $\beta$ and $\gamma$ are $\beta,\gamma \in \{1,i,1+i,i(1+i)\}$.\\
We assert that none of these are possible.\\
Suppose $\beta=1$. Since $\alpha=1$ or $i$, we can write $t=\frac {x^2} {y^2}$ or $t=\frac{ix^2} {y^2}$, i.e. $t^2=\pm \frac{x^4}{y^4}$. Multiplying 
(\ref{eq3}) by $y^4$ we get $x^4\pm y^4=\pm z^2$. It was proved by Hilbert(see \cite{hil}) that this equation has only trivial solutions in Gaussian integers, implying $t=0$.\\
Suppose $\beta=i$. Multiplying (\ref{eq2}) and (\ref{eq3}) we obtain $iy^2=t^3+t$ or $-y^2=t^3+t$, leading to elliptic curves in Weierstrass form $y^2=x^3-x$ and $y^2=x^3+x$ respectively. Using the program \cite{sim}, written in PARI, we compute that the rank of this curve is 0. It is easy to compute, using the Extended Lutz-Nagell Theorem that the torsion subgroup of both these curves over $\mathbb Q(i)$ is $\mathbb Z_2 \oplus \mathbb Z_2$. All the torsion points of these curves satisfy $t(t^4-1)=0$.\\
Suppose $\beta=1+i$ or $i(1+i)$. Multiplying (\ref{eq2}) and (\ref{eq3}) we obtain one of the following three elliptic curves $(1+i)y^2=t^3+t$, $i(1+i)y^2=t^3+t$ and $-(1+i)y^2=t^3+t$. These curves induce the curves $y^2=x^3+2ix$ and $y^2=x^3-2ix$ in Weierstrass form, both of them having rank 0 (again, this is computed using \cite{sim}) and torsion $\mathbb Z_2 \oplus \mathbb Z_2$. All of the torsion points induce $t$ satisfying $t(t^4-1)=0$. Hence, the starting equation has no solutions. \qed

\newtheorem{tm7}[tm]{Lemma}
\begin{tm7} $E(\mathbb Q(i))_{tors}$ cannot be $\mathbb Z_{15}$.
\end{tm7}
\emph{Proof}:\\
By \cite{ra}, case 2.5.4, pages 34 and 35, elliptic curves with torsion subgroups isomorphic to $\mathbb Z_{15}$ are induced by the solutions over $\mathbb Q(i)$ of \begin{equation}s^2+st+s=t^3+t^2\label{eq5}\end{equation} satisfying $$t(t+1)(t^2+t+1)(t^4+3t^3+4t^2+2t+1)(t^4-7t^3-6t^2+2t+1)\neq 0.$$
Using \cite{sim}, we compute that the rank of (\ref{eq5}) over $\mathbb Q(i)$ is 0 and the torsion points give $t=0$ or $-1$, implying that the equation has no solutions.\qed\\
\vspace{0.3cm}\\
\newtheorem{tm8}[tm]{Lemma}
\begin{tm8} $E(\mathbb Q(i))_{tors}$ cannot be $\mathbb Z_2 \oplus \mathbb Z_{12}.$
\end{tm8}
\emph{Proof:}\\
By \cite{ra}, case 2.5.8, pages 42--44, elliptic curves with torsion $\mathbb Z_2 \oplus \mathbb Z_{12}$ are induced by the solutions over $\mathbb Q(i)$ of the equation \begin{equation} s^2=(2t^2-2t+1)(6t^2-6t+1)\label{eq6}\end{equation} satisfying
\begin{equation}t(t-1)(2t-1)(2t^2-2t+1)(3t^2-3t-1)(6t^2-6t+1)\neq 0.\label{uv} \end{equation}
The elliptic curve (\ref{eq6}) has the Weirstrass form
$$y^2=x^3-x^2+x.$$
This curve has rank 0 and $8$ torsion points. They are $\{O,(0,0),(1,\pm 1), (\pm i,\pm 1)\}$ in Weierstrass form, inducing $t=0,1,\frac 1 2$ or $\frac{1\pm i} 2$,  none of them satisfying (\ref{uv}). \qed

\newtheorem{tm9}[tm]{Lemma}
\begin{tm9} $E(\mathbb Q(i))_{tors}$ cannot be $\mathbb Z_{14}$.
\end{tm9}
\emph{Proof:}\\
By \cite{ra}, case 2.5.3, page 31, elliptic curves with torsion $\mathbb Z_{14}$ are induced by the solutions over $\mathbb Q(i)$ of the equation $$s^2+st+s=t^3-t$$ satisfying 
$$t(t^2-1)(t^3-9t^2-t+1)(t^3-2t^2-t+1)\neq 0.$$
The given curve has rank 0 and 6 torsion points over $\mathbb Q(i)$, all of them satisfying $t=0$ or $\pm 1$. \qed\\
\vspace{0.3cm}\\
Lemmas 4, 5, 6, 7 and 8 prove Theorem 1, (ii). Combining this with Lemma 3, we get Theorem 1, (i).
\vspace{0.5cm}\\
\textbf{3 Torsion over $\mathbb Q(\sqrt{-3})$}\\
\\
As some of the proofs in this section are similar to the ones in the previous section, some technical details will be omitted.\\
Let $\omega =\frac{1-\sqrt{-3}}{2}$. It is easy to see that $\omega$ is a primitive sixth root of unity and $\mathbb Z [\omega]$ is the ring of integers of $\mathbb Q(\sqrt{-3})$.
\newtheorem{tm10}[tm]{Lemma}
\begin{tm10} $E(\mathbb Q(\sqrt{-3}))_{tors}$ cannot be $\mathbb Z_{14},\ \mathbb Z_{15},$ or $\mathbb Z_2 \oplus \mathbb Z_{12} $.
\end{tm10}
\emph{Proof:}\\
It is proved in \cite{km} that the torsion subgroup cannot be $\mathbb Z_{14}$(Example 2.5) or $\mathbb Z_2 \oplus \mathbb Z_{12}$(Example 3.2).\\
The proof that the torsion cannot be $\mathbb Z_{15}$ is completely analogous to the proof of Lemma 7.\qed
\newtheorem{tm11}[tm]{Lemma}
\begin{tm11} $E(\mathbb Q(\sqrt{-3}))_{tors}$ cannot be $\mathbb Z_{11}$.
\end{tm11}
\emph{Proof:}\\
As can be seen in \cite{ra}, case 2.5.1, page 25, the solutions $s,t\in\mathbb Q(\sqrt{-3})$ of the equation 
$$s^2-s=t^3-t^2$$
satisfying
$$t(t-1)(t^5-18t^4+35t^3-16t^2-2t+1)\neq 0$$
induce elliptic curves with torsion $\mathbb Z_{11}$ over $ \mathbb Q(\sqrt{-3})$. The rank of this curve is $0$ and there are $5$ torsion points, all of them satisfying $t=0$ or $1$ (see \cite{ra}, Lemma 2.1). \qed
\newtheorem{tm12}[tm]{Lemma}
\begin{tm12} $E(\mathbb Q(\sqrt{-3}))_{tors}$ cannot be $\mathbb Z_2 \oplus \mathbb Z_{10}$.
\end{tm12}
\emph{Proof:}\\
As can be seen in \cite{ra}, case 2.5.7, pages 39--40, the solutions $s,t\in\mathbb Q(\sqrt{-3})$ of the equation 
$$s^2=t^3+t^2-t$$
satisfying
$$t(t^2-1)(t^2-4t-1)(t^2+t-1)\neq 0$$
induce elliptic curves with torsion $\mathbb Z_2 \oplus \mathbb Z_{10}$ over $ \mathbb Q(\sqrt{-3})$. The rank of this curve is $0$ and there are $6$ torsion points, all of them satisfying $t=0,-1$ or $1$ (see \cite{ra}, Lemma 2.4). \qed\\
\vspace{0.3cm}\\
As in the proof of Theorem 1, the hardest part of the proof of Theorem 2 is eliminating the possibility of the torsion being $\mathbb Z_{16}$.
\newtheorem{tm13}[tm]{Lemma}
\begin{tm13} $E(\mathbb Q(\sqrt{-3}))_{tors}$ cannot be $\mathbb Z_{16}$.
\end{tm13}
\emph{Proof:}\\
Again, we have to prove that the equation (\ref{j1}) has no solutions satisfying (\ref{j2}). We follow the same strategy of the proof of Lemma 6, and define $\alpha, \beta$ and $\gamma$ in the same way. It can be proved in the same way as in Lemma 6 that $\alpha$ is a unit and that each of $\beta$ and $\gamma$ is either a unit or twice a unit. As every unit is a square or $\omega$ times a square, and as $\gamma=\alpha\beta \pmod{(\mathbb Q(\sqrt{-3})^*)^2 }$, we see that we have 8 possibilities for the triples $(\alpha,\beta,\gamma)$. If $t$ is a solution of (\ref{j1}), then $t$ has to be the first coordinate of a point on both curves 
$$E_1:\alpha\beta y^2=t^3+t$$ 
and
$$E_2:\alpha\gamma y^2=t^3+2t^2-t.$$
In the following table we give the ranks of these curves depending on $\alpha$ and $\beta$.
\begin{center}
\begin{tabular}{|c|c|c|c|c|}
\hline
$\alpha$ & $\beta$ & $\gamma$ & rank($E_1(\mathbb Q(\sqrt{-3})$) & rank($E_2(\mathbb Q(\sqrt{-3})$)\\
\hline
1 & 1 & 1 & 1 & 0 \\
\hline
1 & $\omega$ & $\omega$ & 1 & 0 \\
\hline
1 & 2 & 2 & 0 & 2 \\
\hline
1 & $2\omega$ & $2\omega$ & 0 & 0 \\
\hline
$\omega$ & $1$ & $\omega$ & 1 & 0 \\
\hline
$\omega$ & $\omega$ & $1$ & 1 & 0 \\
\hline
$\omega$ & $2\omega$ & $2$ & 0 & 0 \\
\hline
$\omega$ & $2$ & $2\omega$ & 0 & 2 \\
\hline
\end{tabular}
\end{center}
As it can be seen, for each case there either $E_1$ or $E_2$ has rank 0, and thus only the torsion points are possible solutions. It remains to check the torsion points of the following four curves:
\begin{equation}
y^2=x^3+2x^2-x,
\label{tor1}
\end{equation}
\begin{equation}
y^2=x^3+2\omega x^2-\omega^2 x,
\label{tor2}
\end{equation}
\begin{equation}
y^2=x^3+4x,
\label{tor3}
\end{equation}
\begin{equation}
y^2=x^3+4\omega^2 x.
\label{tor4}
\end{equation}
The torsion groups were now computed in APECS(\cite{ap}). The torsion of the curves (\ref{tor3}) and (\ref{tor4}) is $\mathbb Z_4$, corresponding to $t=0,\pm 1$ on the curve $E_1$. The torsion of the curves (\ref{tor1}) and (\ref{tor2}) is $\mathbb Z_2$, corresponing to $t=0$ on the curve $E_2$.\\
We obtain that none of the $t$ induced by the torsion points satisfies (\ref{j2}). We conclude that the torsion $\mathbb Z_{16}$ is impossible.\qed\\
\vspace{0.3cm}\\
Lemmas 9, 10, 11 and 12 combined with Lemma 3 prove Theorem 2.\\
\vspace{0.2 cm}\\
\textbf{Remark}\\
In order to prove that there are no elliptic curves with torsion $\mathbb Z_{13}$ over $\mathbb Q(i)$ and $\mathbb Q(\sqrt{-3})$, one would have to prove that there are no solutions in the respective quadratic field of the equation \begin{equation}
s^2=t^6-2t^5+t^4-2t^3+6t^2-4t+1
\label{he1}
\end{equation}
satisfying
$$t(t-1)(t^3-4t^2+t+1)\neq 0.$$
Similarly, to prove that there are no elliptic curves with torsion $\mathbb Z_{18}$ over $\mathbb Q(\sqrt{-3})$ one would have to prove that there do not exist $s,t\in \mathbb Q(\sqrt{-3})$ satisfying
\begin{equation}
s^2=t^6+2t^5+5t^4+10t^3+10t^2+4t+1
\label{he2}
\end{equation}
and
$$t(t+1)(t^2+t+1)(t^3-3t-1)\neq 0.$$
Note that (\ref{he1}) and (\ref{he2}) are both hyperelliptic curves of genus 2.

\vspace{0.2 cm}

\textbf{Acknowledgements}\\
The author would like to thank Andrej Dujella for motivating and helpful discussions. Also, the author is grateful to Mirela Juki\'c-Bokun and the referee for pointing out some mistakes in earlier versions of this manuscript.

\small{FILIP NAJMAN}\\
\small{DEPARTMENT OF MATHEMTICS,\\ UNIVERSITY OF ZAGREB,\\ BIJENI\v CKA CESTA 30, 10000 ZAGREB,\\ CROATIA}\\
\emph{E-mail address:} fnajman@math.hr

\begin{thebibliography}{1}
\bibitem{ap}
I. Connell, \emph{APECS}, ftp://ftp.math.mcgill.ca/pub/apecs/.


\bibitem{fu2}
Y.  Fujita,  \emph{Torsion subgroups of elliptic curves in elementary
abelian $2$-extensions of $\mathbb Q$},
J. Number Theory  \textbf{114}  (2005),  124--134.
\bibitem{hil}
D. Hilbert,
\emph{Jahresbericht d. Deutschen Math.-Vereinigung}, 4, 1894/1895, 517--525.

\bibitem{jkp}
D. Jeon, C. H. Kim and E. Park,
\emph{On the torsion of elliptic curves over quartic number fields},
J. London Math. Soc. (2) \textbf{74} (2006), 1--12.
\bibitem{kam}
S. Kamienny, \emph{Torsion points on elliptic curves and $q$-coefficients of modular forms}, Invent. Math. \textbf{109} (1992), 221--229.
\bibitem{km}
M. A. Kenku and F. Momose, \emph{Torsion points on  elliptic curves defined over quadratic
fields}, Nagoya Math. J. \textbf{109} (1988), 125--149.
\bibitem{kn}
A. Knapp,
\emph{Elliptic curves}, Princeton Univ. Press, 1992.

\bibitem{ll}
M. Laska and M. Lorenz, 
\emph{Rational points on elliptic curves over $\mathbb Q$
in elementary abelian $2$-extensions of $\mathbb Q$},  J. Reine Angew. Math. \textbf{355}
(1985), 163--172.
\bibitem{ra}
F. P. Rabarison,
\emph{Torsion et rang des courbes elliptiques definies sur les corps de nombres algebriques}, Doctorat de Universite de Caen, 2008.

\bibitem{re}
M. A. Reichert,
\emph{Explicit Determination of Nontrivial Torsion Structures of Elliptic Curves Over Quadratic
Number Fields}, 
Math. Comp. \textbf{174} (1986), 637--658.

\bibitem{sim}  
D. Simon,
\emph{Le fichier gp},
http://www.math.unicaen.fr/$\sim$simon/ell.gp.
\bibitem{th}
T. Thongjunthug,
\emph{Elliptic curves over $\mathbb Q(i)$,} Honours thesis (2006).
\end{thebibliography}
\end{document}